\documentclass[11pt]{article}
\usepackage[margin=1in]{geometry}

\usepackage{amsmath} %
\usepackage{amssymb}  %
\usepackage{amsthm}
\usepackage{color}

\theoremstyle{plain}  %
\newtheorem{theorem}{Theorem}[section]
\newtheorem{lemma}[theorem]{Lemma}
\newtheorem{proposition}[theorem]{Proposition}
\newtheorem{corollary}[theorem]{Corollary}
\newtheorem{definition}[theorem]{Definition}

\theoremstyle{definition}

\theoremstyle{remark} %

\newtheorem{remark}{Remark}[section]

\def\0{{\bf 0}}
\def\1{{\bf 1}}
\def\alex#1{{#1}}

\long\def\jnt#1{{#1}}
\long\def\old#1{}
\long\def\jntn#1{{#1}}

\definecolor{DarkerGreen}{RGB}{0,170,0}
\long\def\aaa#1{{#1}}
\definecolor{orange}{rgb}{1,0.5,0}
\long\def\aaan#1{{#1}}

\title{\LARGE \bf NP-hardness of Deciding Convexity of \\Quartic Polynomials and Related Problems}
 \author{Amir Ali Ahmadi, Alex Olshevsky, Pablo A. Parrilo, and John N. Tsitsiklis \thanks{The authors are with the Laboratory for Information and
Decision Systems, Department of Electrical Engineering and
Computer Science, Massachusetts Institute of Technology. Email:
\{\texttt{a\_a\_a}, \texttt{alex\_o}, \texttt{parrilo},
\texttt{jnt}\}\texttt{@mit.edu}.  } \thanks{\aaa{This research was
partially supported by the NSF Focused Research Group Grant on
Semidefinite Optimization and Convex Algebraic Geometry
DMS-0757207 and by the NSF grant ECCS-0701623.}} 
}
\begin{document}
\date{}
\maketitle

\begin{abstract}
\noindent We show that unless P=NP, there exists no polynomial
time (or even pseudo-polynomial time) algorithm that can decide
whether a multivariate polynomial of degree four (or higher even
degree) is \aaa{globally} convex. This solves a problem that has
been open since 1992 when \aaa{N. Z.} Shor \alex{asked} for the
complexity of deciding convexity for quartic polynomials. We also
prove that deciding strict convexity, strong convexity,
quasiconvexity, and pseudoconvexity of polynomials of even degree
four or higher is strongly NP-hard. By contrast, we show that
quasiconvexity and pseudoconvexity of odd degree polynomials can
be decided in polynomial time.
\end{abstract}

\section{Introduction}
The role of \emph{convexity} in modern day mathematical
programming has proven to be remarkably fundamental, to the point
that tractability of an optimization problem is nowadays assessed,
more often than not, by whether or not the problem benefits from
some sort of underlying convexity. In the famous words of
Rockafellar~\cite{Roc_SIAM_Lagrange}: \vspace{-5pt}
\begin{itemize}
\item[]
``In fact the great watershed in optimization isn't between
linearity and nonlinearity,\\ but convexity and nonconvexity.''
\end{itemize}
\vspace{-5pt}
But how easy is it to distinguish between convexity and
nonconvexity? Can we decide in an efficient manner if a given
optimization problem is convex?

A class of optimization problems that allow for a rigorous study
of this question from a computational complexity viewpoint is the
class of polynomial optimization problems. \alex{These are}
optimization problems where the objective is given by a polynomial
function and the feasible set is described by polynomial
inequalities. Our research in this direction was motivated by a
concrete question of N. Z. Shor that appeared as one of seven open
problems in complexity theory for numerical optimization put
together by Pardalos and Vavasis in 1992~\cite{open_complexity}:
\vspace{-5pt}
\begin{itemize}
\item[]
``Given a degree-$4$ polynomial in $n$ variables, what is the
complexity of determining\\ whether this polynomial describes a
convex function?''
\end{itemize}
\vspace{-5pt}
As we will explain in more detail shortly, the reason why Shor's
question is specifically about degree $4$ polynomials is that
deciding convexity of odd degree polynomials is trivial and
deciding convexity of degree $2$ (quadratic) polynomials can be
reduced to the simple task of checking \jnt{whether} a constant matrix is
positive semidefinite. So, the first interesting case really
occurs for degree $4$ (quartic) polynomials. Our main contribution
in this paper (Theorem~\ref{thm:convexity.quartic.nphard} in
Section~\ref{subsec:convexity.hard.degrees}) is to show that
deciding convexity of polynomials is strongly NP-hard already for
polynomials of degree $4$.

The implication of NP-hardness of this problem is that
\aaan{unless P=NP}, there exists no algorithm that can take as
input the (\alex{rational}) coefficients of a quartic polynomial,
have running time bounded by a polynomial in the number of bits
needed to represent the coefficients, and output correctly on
every instance whether or not the polynomial is convex.
Furthermore, the fact that our NP-hardness result is in the strong
sense (as opposed to weakly NP-hard problems such as KNAPSACK)
implies, roughly speaking, that the problem remains NP-hard even
when the magnitude of the coefficients of the polynomial are
restricted to be ``small.'' For a strongly NP-hard problem, even a
pseudo-polynomial time algorithm cannot exist unless P=NP.
See~\cite{GareyJohnson_Book} for precise definitions and more
details.

\aaan{There are many areas of application where one would like to
establish convexity of polynomials.} Perhaps the simplest example
is in global minimization of polynomials, where it could be very
useful to decide first \jnt{whether the polynomial to be optimized
is convex. Once convexity is verified,} then every local minimum
is global and very basic techniques (e.g., gradient descent) can
find a global minimum---a task that is in general NP-hard in
\aaan{the} absence of
convexity~\cite{Minimize_poly_Pablo},~\cite{nonnegativity_NP_hard}.
As another example, if we can certify that a homogeneous
polynomial is convex, \aaa{then we define a gauge (or Minkowski)
norm based on its convex sublevel sets, which may be useful in
many applications}. In several other problems \jnt{of practical
relevance,} we might not just be interested in checking whether a
given polynomial is convex, but to \emph{parameterize} a family of
convex polynomials and perhaps search or optimize over them. For
example we might be interested in approximating the convex
envelope of a complicated nonconvex function with a convex
polynomial, or in fitting a convex polynomial to a set of data
points with minimum error~\cite{convex_fitting}. Not surprisingly,
if testing membership to the set of convex polynomials is hard,
searching and optimizing over \jnt{that set} also turns out to be
a hard problem.

We also extend our hardness result to some variants of convexity,
namely, the problems of deciding \emph{strict convexity},
\emph{strong convexity}, \emph{pseudoconvexity}, and
\emph{quasiconvexity} of polynomials. Strict convexity is a
property that is often useful to check because it guarantees
uniqueness of the optimal solution in optimization problems. The
notion of strong convexity is a common assumption in convergence
analysis of many iterative Newton-type algorithms in optimization
theory; see, e.g.,~\cite[Chaps.\ 9--11]{BoydBook}. So, in order to
ensure the theoretical convergence rates promised by many of these
algorithms, one needs to first make sure that the objective
function is strongly convex. The problem of checking
quasiconvexity (convexity of sublevel sets) of polynomials also
arises frequently in practice. For instance, if the feasible set
of an optimization problem is defined by polynomial inequalities,
by certifying quasiconvexity of the defining polynomials we can
ensure that the feasible set is convex. In several statistics and
clustering problems, we are interested in finding minimum volume
convex sets that contain a set of data points in space. This
problem can be tackled by searching over the set of quasiconvex
polynomials~\cite{convex_fitting}. In economics, quasiconcave
functions are prevalent as desirable utility
functions~\cite{Testing_convexity_economics},~\cite{Quasiconcave_programming}.
In control and systems theory, it is useful at times to search for
quasiconvex Lyapunov functions whose convex sublevel sets contain
relevant information about the trajectories of a dynamical
system~\cite{Chesi_Hung_journal},~\cite{AAA_PP_CDC10_algeb_convex}.
Finally, the notion of pseudoconvexity is a natural generalization
of convexity that inherits many of the attractive properties of
convex functions. For example, every stationary point or every
local minimum of a pseudoconvex function must be a global minimum.
Because of these nice features, pseudoconvex programs have been
studied extensively in nonlinear
programming~\cite{Mangasarian_Pseudoconvex_fns},%
~\cite{pseudoconvex_nonnegative_vars}.

As an outcome of close to a century of research in convex
analysis, numerous necessary, sufficient, and exact conditions for
convexity and all of its variants are available; see,
e.g.,~\cite[Chap.~3]{BoydBook},~\cite{Second_order_pseudoconvexity},~\cite{Matrix_theoretic_quasiconvexity},~\cite{Criteria_comparison_quasiconvexity_pseudoconvexity},~\cite{Testing_convexity_economics},
~\cite{NLP_Book_Mangasarian} and references therein for a by no
means exhaustive list. Our results %
suggest that none
of the exact characterizations of these notions can be efficiently
checked for polynomials. In fact, when turned upside down, many of
these equivalent formulations reveal new NP-hard problems; see, e.g.,
Corollary~\ref{cor:nonnegativity.quartic.nphard} and
\ref{cor:semialgeb.set.nphard}.

\subsection{Related Literature}
There are several results in the literature \jnt{on} the complexity of
various special cases of polynomial optimization problems. The
interested reader can find many of these results in the edited volume
of Pardalos~\cite{Pardalos_Book_Complexity_num_opt} or in the survey
papers of de Klerk~\cite{deKlerk_complexity_survey}, and Blondel and
Tsitsiklis~\cite{BlTi1}. A very general and fundamental concept in
certifying feasibility of polynomial equations and inequalities is the
Tarski--Seidenberg quantifier elimination
theory~\cite{Tarski_quantifier_elim},~\cite{Seidenberg_quantifier_elim},
from which it follows that all of the problems that we consider in
this paper are algorithmically \emph{decidable}. This means that there
are algorithms that on all instances of our problems of interest halt
in finite time and always output the correct yes--no
answer. Unfortunately, algorithms based on quantifier elimination or
similar decision algebra techniques have running times that are at
least exponential in the number of variables
\cite{Algo_real_algeb_geom_Book}, and in practice can only solve
problems with very few parameters.

When we turn to the issue of polynomial time solvability, perhaps
the most relevant result for our purposes is the NP-hardness of
deciding nonnegativity of quartic polynomials and biquadratic
forms (see Definition~\ref{defn:biquad.forms}); the reduction that
we give in this paper will in fact be from the latter problem.
\jnt{As we will see in
Section~\ref{subsec:convexity.hard.degrees}, it turns out that
deciding convexity of quartic forms is \aaan{equivalent to}
checking nonnegativity of a special \aaan{class} of biquadratic
forms, which are themselves a special \aaan{class} of quartic
forms.} The NP-hardness of checking nonnegativity of quartic forms
follows, e.g., as a direct consequence of NP-hardness of testing
matrix copositivity, a result proven by Murty and
\aaa{Kabadi}~\cite{nonnegativity_NP_hard}. As for the hardness of
checking nonnegativity of biquadratic forms, we know of two
different proofs. The first one is due to
Gurvits~\cite{Gurvits_quantum_entag_hard}, \jnt{who} proves that
the entanglement problem in quantum mechanics (i.e., the problem
of distinguishing separable quantum states from entangled ones) is
NP-hard. A dual \jnt{reformulation of this} %
result shows directly that checking nonnegativity of biquadratic
forms is NP-hard; see~\cite{Pablo_Sep_Entang_States}. The second
proof is due to Ling et al.~\cite{Ling_et_al_Biquadratic},
\jnt{who} use a theorem of Motzkin and Straus to give a very short
and elegant reduction from the maximum independent set problem in
graphs.

The only work in the literature on the hardness of deciding
polynomial convexity that we are aware of is the work of Guo on the
complexity of deciding convexity of quartic polynomials over
simplices~\cite{complexity_simplex_convexity}. Guo discusses some of
the difficulties that arise from this problem, but he does not prove that deciding convexity
of polynomials over simplices is NP-hard.
Canny shows
in~\cite{Canny_PSPACE} that the existential theory of the real
numbers can be decided in PSPACE. From this, it follows that testing
several properties of polynomials, including nonnegativity and
convexity, can be done in polynomial space.
In~\cite{Nie_PMI_SDP_repres.}, Nie proves that the related notion of
\emph{matrix convexity} is NP-hard for polynomial matrices whose
entries are quadratic forms. %

On the \alex{algorithmic} side, several techniques have been
proposed both for testing convexity of sets and convexity of
functions. \aaan{Rademacher and Vempala present and analyze
randomized algorithms for testing the relaxed notion of
\emph{approximate convexity}~\cite{Test_Geom_Convexity}. 
In~\cite{Lasserre_set_convexity}, Lasserre proposes a semidefinite
programming hierarchy for testing convexity of basic closed
semialgebraic sets; a problem that we also prove to be NP-hard
(see Corollary~\ref{cor:semialgeb.set.nphard}).}
As for testing convexity of functions, an approach that some
convex optimization parsers (e.g.,  \texttt{CVX}~\cite{cvx}) take
is to start with some ground set of convex functions and then
check \aaan{whether} the desired function can be obtained
\aaan{by} applying a set of convexity preserving operations to the
functions in the ground
set~\cite{Crusius_thesis},~\cite[p.~79]{BoydBook}.
\aaan{Techniques of this type that are based on the calculus of
convex functions are successful for a large range of applications.
However, when applied to general polynomial functions, they can
only detect a subclass of convex polynomials.}

\jnt{Related} to convexity of polynomials, a concept that has
attracted recent attention is the algebraic notion of
\emph{sos-convexity} (see
Definition~\ref{defn:sos-convex})~\cite{Helton_Nie_SDP_repres_2},~\cite{Lasserre_Jensen_inequality},~\cite{Lasserre_Convex_Positive},~\cite{AAA_PP_CDC10_algeb_convex},~\cite{convex_fitting},~\cite{Chesi_Hung_journal}.
This is a powerful sufficient condition for convexity \aaa{that}
relies on an appropriately defined sum of squares decomposition of
the Hessian matrix, and can be efficiently checked by solving a
single semidefinite program. However,
in~\cite{AAA_PP_not_sos_convex_journal},~\cite{AAA_PP_CDC09_HessianNotFactor},
Ahmadi and Parrilo gave an explicit counterexample to show that
not every convex polynomial is sos-convex. \aaa{The NP-hardness
result in this work certainly justifies the existence of such a
counterexample and more generally suggests that \emph{any}
polynomial time algorithm attempted for checking polynomial
convexity is doomed to fail on some hard instances.}

\subsection{\alex{Contributions and} organization of the paper}
\jnt{The} \alex{main} contribution of this paper is to
\aaa{establish} the computational complexity of deciding
convexity, strict convexity, strong convexity, pseudoconvexity,
and quasiconvexity of polynomials for any given degree. (See
Table~1 in Section~\ref{sec:summary.conclusions} for a quick
summary.) The results are mainly divided in three sections, with
Section~\ref{sec:convexity} covering convexity,
Section~\ref{sec:strict.strong} covering strict and strong
convexity, and Section~\ref{sec:quasi.pseudo} covering
quasiconvexity and pseudoconvexity. \aaa{These three sections
follow a similar \aaan{pattern} and are each divided into three
parts: first, the definitions and basics, second, the degrees for
which the questions can be answered in polynomial time, and third,
the degrees for which the questions are NP-hard.}

Our main reduction, which establishes NP-hardness of checking
convexity of quartic forms, is given in
Section~\ref{subsec:convexity.hard.degrees}. This hardness result
is extended to strict and strong convexity in
Section~\ref{subsec:strict.strong.hard.degrees}, and to
quasiconvexity and pseudoconvexity in
Section~\ref{subsec:quasi.pseudo.hard.degrees}. \aaa{By contrast,}
we show \aaa{in Section~\ref{subsec:quasi.pseudo.easy.degrees}}
that quasiconvexity and pseudoconvexity of odd degree polynomials
can be decided in polynomial time. Finally, a summary of our
results and some concluding remarks are presented in
Section~\ref{sec:summary.conclusions}.

\section{Complexity of deciding convexity}\label{sec:convexity}
\subsection{Definitions and basics}\label{subsec:convexity.basics}
A \aaa{(multivariate)} \emph{polynomial} $p(x)$ in variables
$x\mathrel{\mathop:}=(x_1,\ldots,x_n)^T$ is a function from
$\mathbb{R}^n$ to $\mathbb{R}$ that is a finite linear combination
of monomials:
\begin{equation}
p(x)=\sum_{\alpha}c_\alpha x^\alpha=\sum_{\alex{\alpha_1, \ldots,
\alpha_n}} c_{\alex{\alpha_1,\ldots,\alpha_n}} x_1^{\alpha_1} \cdots
x_n^{\alpha_n} ,
\end{equation}
\alex{where the sum is over \aaa{$n$-tuples of} nonnegative
integers $\alpha_i$}. An algorithm for testing some property of
polynomials will have as its input an ordered list of the
coefficients $c_\alpha$. Since our complexity results are based on
models of digital computation, where the input must be represented
by a finite number of bits, the coefficients $c_\alpha$ for us
will always be rational numbers, which upon clearing the
denominators can be taken to be integers. So, for the remainder of
the paper, even when not explicitly stated, we will always have
$c_\alpha \in \mathbb{Z}$.

The \emph{degree} of a monomial $x^\alpha$ is equal to $\alpha_1 +
\cdots + \alpha_n$. The degree of a polynomial $p(x)$ is defined to
be the highest degree of its component monomials. A simple counting
argument shows that a polynomial of degree $d$ in $n$ variables has
$\binom{n+d}{d}$ coefficients. A \emph{homogeneous polynomial} (or a
\emph{form}) is a polynomial where all the monomials have the same
degree. A form $p(x)$ of degree $d$ is a homogeneous function of
degree $d$ (since it satisfies $p(\lambda x)=\lambda^d p(x)$), and
has $\binom{n+d-1}{d}$ coefficients.

A polynomial $p(x)$ is said to be \emph{nonnegative} or
\emph{positive semidefinite (psd)} if $p(x)\geq0$ for all
$x\in\mathbb{R}^n$. Clearly, a necessary condition for a polynomial
to be psd is for its degree to be even. We say that $p(x)$ is a
\emph{sum of squares (sos)}, if there exist polynomials
$q_{1}(x),\ldots,q_{m}(x)$ such that
$p(x)=\sum_{i=1}^{m}q_{i}^{2}(x)$. Every sos polynomial is obviously
psd. A \emph{polynomial matrix} $P(x)$ is a matrix with polynomial
entries. We say that a polynomial matrix $P(x)$ is \emph{PSD}
(denoted $P(x)\succeq0$) if it is positive semidefinite in the
matrix sense for every value of the indeterminates $x$.
\jnt{(Note the upper case convention for matrices.)} It is easy
to see that $P(x)$ is PSD if and only if the scalar polynomial
$y^TP(x)y$ in variables $(x;y)$ is psd.

\aaa{We recall that a polynomial $p(x)$ is convex if and only if
its Hessian matrix, which will be generally denoted by $H(x)$, is
PSD.}

\subsection{Degrees that are easy}\label{subsec:convexity.easy.degrees}
The question of deciding convexity is trivial \aaa{for odd degree
polynomials}. Indeed, it is easy to check that linear polynomials
($d=1$) are always convex and that polynomials of odd degree
$d\geq3$ can never be convex. \aaa{The case of quadratic
polynomials ($d=2$) is also straightforward.}
\aaa{A quadratic polynomial $p(x)=\frac{1}{2}x^TQx+q^Tx+c$ is
convex if and only if the constant matrix $Q$ is positive
\jnt{semidefinite}. This can be decided in polynomial time for
example by performing Gaussian pivot steps along the main diagonal
of $Q$~\cite{nonnegativity_NP_hard} or by computing the
characteristic polynomial of $Q$ exactly and then checking that
the signs of its coefficients alternate~\cite[p.
403]{HJ_Matrix_Analysis_Book}.}

Unfortunately, the results that come next suggest that the case of
quadratic polynomials is essentially the only nontrivial case
where convexity can be efficiently \aaa{decided}.

\subsection{Degrees that are hard}\label{subsec:convexity.hard.degrees}
The main hardness result of the paper is the following theorem.
\begin{theorem}\label{thm:convexity.quartic.nphard}
Deciding convexity of degree four polynomials is \jnt{strongly} NP-hard. This is
true even when the polynomials are restricted to be homogeneous.
\end{theorem}
We will give a reduction from the problem of deciding
nonnegativity of biquadratic forms. We start by recalling some
basic facts about biquadratic forms and sketching the idea of the
proof.

\begin{definition}\label{defn:biquad.forms}
A \emph{biquadratic form} \alex{$b(x;y)$ is a \aaa{form in the
variables $x=(x_1, \ldots, x_n)^T$ and $y=(y_1, \ldots, y_m)^T$
that} can be written as}
\begin{equation}\label{eq:defn.biquad.form}
b(x;y)\old{\mathrel{\mathop:}}=\sum_{i\leq j, \, k\leq
l}\alpha_{ijkl}x_ix_jy_ky_l.
\end{equation}
\end{definition}
\aaa{Note that for fixed $x$, $b(x;y)$ becomes a quadratic form in
$y$, and for fixed $y$, it becomes a quadratic form in $x$.} Every
biquadratic form is a quartic form, but the converse is of course
not true. It follows from a result of Ling et
al.~\cite{Ling_et_al_Biquadratic} that deciding nonnegativity of
biquadratic forms is strongly NP-hard even when $n=m$. In the
sequel, we will always have $n=m$.

It is not difficult to see that any biquadratic form $b(x;y)$ can
be written in the form \begin{equation} \label{eq:b_represent}
b(x;y)=y^TA(x)y\end{equation} (or of course as $x^TB(y)x$) for
some symmetric polynomial matrix $A(x)$ whose entries are
quadratic forms. Therefore, it is strongly NP-hard to decide
whether a symmetric polynomial matrix with quadratic form entries
is PSD. \alex{One might hope that this would lead to a quick proof
\aaa{of NP-hardness of testing convexity of quartic forms},
because} the Hessian of a quartic form is exactly a symmetric
polynomial matrix with quadratic form entries. \aaa{However, the
major problem that stands in the way is that not every polynomial
matrix is a \emph{valid Hessian}. Indeed, if any of the partial
derivatives between the entries of $A(x)$ do not commute (e.g., if
$\frac{\partial{A_{11}(x)}}{\partial{x_2}}\neq\frac{\partial{A_{12}(x)}}{\partial{x_1}}$),
then $A(x)$ cannot be the matrix of second derivatives of some
polynomial. This is because all mixed third partial derivatives of
polynomials must commute.}

Our task is therefore to prove that even with these additional
constraints \jnt{on} the entries of $A(x)$, the problem of
deciding positive semidefiniteness of such matrices remains
NP-hard. \aaa{We will show} that any given symmetric $n \times n$
matrix $A(x)$, whose entries are quadratic forms, \jnt{can be
embedded} in a $2n \times 2n$ polynomial matrix $H(x,y)$, again
with quadratic form entries, \jnt{so} that $H(x,y)$ is a valid
Hessian and $A(x)$ is PSD if and only if $H(x,y)$ is. In fact,
\alex{we will directly construct the polynomial $f(x,y)$ whose
Hessian is the matrix $H(x,y)$.} \aaa{This is done in the next
theorem, which establishes the correctness of our main reduction.
Once this theorem is proven, the proof of
Theorem~\ref{thm:convexity.quartic.nphard} will become immediate.}

\aaa{
\begin{theorem}\label{thm:main.reduction}
 \alex{Given a biquadratic form $b(x;y)$,
define the} the $n
\times n$ polynomial matrix $C(x,y)$ by setting
\begin{equation}\label{eq:C(x,y).defn.}
\alex{[}C(x,y)\alex{]_{ij}}\mathrel{\mathop:}=\frac{\partial{b(x;y)}}{\partial{x_{\alex{i}}}\partial{y_{\alex{j}}}},
\end{equation} and let $\gamma$ be the largest coefficient, in
absolute value, \alex{of any monomial present in some entry of the matrix
$C(x,y)$.}
Let $f$ be the form given by
\begin{equation}\label{eq:f(x,y).defn.}
f(x,y)\mathrel{\mathop:}=b(x;y) + \frac{n^2\gamma}{2}\Big(
\sum_{i=1}^n x_i^4+\sum_{i=1}^n
y_i^4+\sum_{\substack{i,j=1,\ldots,n\\i\aaan{<} j}}
x_i^2x_j^2+\sum_{\substack{i,j=1,\ldots,n\\i\aaan{<} j}}
y_i^2y_j^2\Big).\end{equation}
Then, $b(x;y)$ is psd if and only if $f(x,y)$ is convex.
\end{theorem}
}

\begin{proof}

\aaa{Before we prove the claim, let us make a few observations and
try to shed light on the intuition behind this construction}. We
will use $H(x,y)$ to denote the Hessian of $f$. This is a $2n
\times 2n$ polynomial matrix \alex{whose entries are quadratic
forms.} The polynomial $f$ is convex if and only if $z^TH(x,y)z$
is psd. For bookkeeping purposes, let us split the variables $z$
as $z\mathrel{\mathop:}=(z_x, z_y)^T$, where $z_x$ and $z_y$ each
belong to $\mathbb{R}^n$. It will also be helpful to give a name
to the second group of terms in the definition of $f(x,y)$ in
(\ref{eq:f(x,y).defn.}). So, let
\begin{equation}\label{eq:g(x,y).defn.}
g(x,y)\mathrel{\mathop:}=\frac{n^2\gamma}{2}\aaa{\Big(
\sum_{i=1}^n x_i^4+\sum_{i=1}^n
y_i^4+\sum_{\substack{i,j=1,\ldots,n\\i\aaan{<} j}}
x_i^2x_j^2+\sum_{\substack{i,j=1,\ldots,n\\i\aaan{<} j}}
y_i^2y_j^2\Big)}.
\end{equation}
We denote the Hessian matrices of $b(x,y)$ and $g(x,y)$ with
$H_b(x,y)$ and $H_g(x,y)$ respectively. Thus,
$H(x,y)=H_b(x,y)+H_g(x,y)$. \aaa{Let us first focus on the
structure of $H_b(x,y)$}. \alex{Observe that if we define}
$$\alex{[}A(x)\alex{]_{ij}}=\frac{\partial{b(x;y)}}{\partial{y_{\alex{i}}}\partial{y_{\alex{j}}}},$$
\alex{then $A(x)$ depends only on $x$, and  \begin{equation}
 \frac{1}{2} y^T A(x) y = b(x;y). \label{eq:yAy=b} \end{equation} }
\alex{Similarly, if we let}
$$\alex{[}B(y)\alex{]_{ij}}=\frac{\partial{b(x;y)}}{\partial{x_{\alex{i}}}\partial{x_{\alex{j}}}},$$
\alex{then $B(y)$ depends only on $y$, and \begin{equation}\label{eq:x.B.x=b}
\frac{1}{2}x^TB(y)x=b(x;y).
\end{equation}}\alex{From Eq.\ (\ref{eq:yAy=b}), we have that $b(x;y)$ is psd if and only if $A(x)$ is \aaa{PSD}; from  Eq.\ (\ref{eq:x.B.x=b}),
we see that $b(x;y)$ is psd if and only if $B(y)$ is \aaa{PSD}. }

Putting the blocks together, we have
\begin{equation}\label{eq:Hb}
H_b(x,y)=\begin{bmatrix} B(y) & C(x,y) \\ C^T(x,y) & A(x)
\end{bmatrix}.
\end{equation} The matrix $C(x,y)$ is not in general symmetric.
\aaan{The entries of $C(x,y)$ consist of square-free monomials
that are each a multiple of} $x_iy_j$ for some \jnt{$i$, $j$,
with} $1\leq i,j \leq n$; (see (\ref{eq:defn.biquad.form}) and
(\ref{eq:C(x,y).defn.})).

\aaa{The Hessian $H_g(x,y)$ of the polynomial $g(x,y)$ in
(\ref{eq:g(x,y).defn.}) is given by}

\begin{equation}\label{eq:Hg}
H_g(x,y)=\frac{n^2\gamma}{2}\begin{bmatrix}H_g^{11}(x) & 0 \\ 0 &
H_g^{22}(y)
\end{bmatrix},
\end{equation}
where
\begin{equation}
H_g^{11}(x)=\begin{bmatrix}
12x_1^2+2\aaa{\displaystyle\sum_{\substack{i=1,\ldots,n\\i\neq 1}}} x_i^2 & 4x_1x_2 & \cdots & 4x_1x_n  \\
4x_1x_2 & 12x_2^2+2\aaa{\displaystyle\sum_{\substack{i=1,\ldots,n\\i\neq 2}}} x_i^2  & \cdots & 4x_2x_n \\
\vdots & \aaan{\vdots} &\ddots &\vdots \\
4x_1x_n &  \cdots &4x_{n-1}x_n & 12x_n^2+2\aaa{\displaystyle\sum_{\substack{i=1,\ldots,n\\i\neq n}}} x_i^2  \\
\end{bmatrix},
\end{equation} and \begin{equation}
H_g^{22}(y)=\begin{bmatrix}
12y_1^2+2\aaa{\displaystyle\sum_{\substack{i=1,\ldots,n\\i\neq 1}}} y_i^2 & 4y_1y_2 & \cdots & 4y_1y_n  \\
4y_1y_2 & 12y_2^2+2\aaa{\displaystyle\sum_{\substack{i=1,\ldots,n\\i\neq 2}}} y_i^2  & \cdots & 4y_2y_n \\
\vdots & \aaan{\vdots} &\ddots &\vdots \\
4y_1y_n &  \cdots &4y_{n-1}y_n & 12y_n^2+2\aaa{\displaystyle\sum_{\substack{i=1,\ldots,n\\i\neq n}}} y_i^2  \\
\end{bmatrix}.
\end{equation}
\aaa{Note that all diagonal elements of $H_g^{11}(x)$ \aaan{and}
$H_g^{22}(y)$ contain the square of every variable $x_1, \ldots,
x_n$ \aaan{and} $y_1, \ldots, y_n$ \aaan{respectively}.}

\aaa{\aaan{We fist give an intuitive summary of the rest of the
proof.} If $b(x;y)$ is not psd, then $B(y)$ and $A(x)$ are not PSD
and hence $H_b(x,y)$ is not PSD. Moreover, adding $H_g(x,y)$ to
$H_b(x,y)$ cannot help make $H(x,y)$ PSD because the dependence of
the diagonal blocks of $H_b(x,y)$ and $H_g(x,y)$ on $x$ and $y$
\aaan{runs} backwards. On the other hand, if $b(x;y)$ is psd, then
$H_b(x,y)$ will have PSD diagonal blocks. \aaan{In principle,
$H_b(x,y)$ might} still not be PSD because of the off-diagonal
block $C(x,y)$. However, the squares in the diagonal elements of
$H_g(x,y)$ \aaan{will be shown to} dominate the monomials of
$C(x,y)$ \aaan{and} make $H(x,y)$ PSD.}

Let us now prove the theorem formally. One direction is easy: if
$b(x;y)$ is not psd, then $f(x,y)$ is not convex. Indeed, if there
exist $\bar{x}$ and $\bar{y}$ in $\mathbb{R}^n$ such that
$b(\bar{x};\bar{y})<0$, then
$$z^TH(x,y)z \Big \vert_{z_x=0, x=\bar{x}, y=0,
z_y=\bar{y}}=\bar{y}^TA(\bar{x})\bar{y}=2b(\bar{x};\bar{y})<0.$$

\aaan{For the converse,} suppose that $b(x;y)$ is psd; we will
prove that $z^TH(x,y)z$ is psd and hence $f(x,y)$ is convex. We
have
\begin{equation}
z^TH(x,y)z=z_y^TA(x)z_y+z_x^TB(y)z_x+2z_x^TC(x,y)z_y+\frac{n^2\gamma}{2}z_x^TH_g^{11}(x)z_x+\frac{n^2\gamma}{2}z_y^TH_g^{22}(y)z_y.
\end{equation}
Because $z_y^TA(x)z_y$ and $z_x^TB(y)z_x$ are psd by assumption
(see \jnt{(\ref{eq:yAy=b}) and} (\ref{eq:x.B.x=b})), it suffices to show that
$z^TH(x,y)z-z_y^TA(x)z_y-z_x^TB(y)z_x$ is psd. In fact, we will
show that $z^TH(x,y)z-z_y^TA(x)z_y-z_x^TB(y)z_x$ is a sum of
squares.

After \jnt{some} regrouping of terms we can write
\begin{equation}\label{eq:p1+p2+p3}
z^TH(x,y)z-z_y^TA(x)z_y-z_x^TB(y)z_x=p_1(x,y,z)+p_2(x,z_x)+p_3(y,z_y),
\end{equation}
where
\begin{equation}\label{eq:p1}
p_1(x,y,z)=2z_x^TC(x,y)z_y+n^2\gamma
\Big(\sum_{\aaa{i=1}}^nz_{x,i}^2\Big)\Big(\sum_{\aaa{i=1}}^nx_i^2\Big)+n^2\gamma\Big(\sum_{\aaa{i=1}}^nz_{y,i}^2\Big)\Big(\sum_{\aaa{i=1}}^ny_i^2\Big),
\end{equation}
\begin{equation}\label{eq:p2}
p_2(x,z_x)=n^2\gamma z_x^T\begin{bmatrix}
5x_1^2& 2x_1x_2 & \cdots & 2x_1x_n  \\
2x_1x_2 & 5x_2^2  & \cdots & 2x_2x_n \\
\vdots & \aaan{\vdots} &\ddots &\vdots \\
2x_1x_n &  \cdots &2x_{n-1}x_n & 5x_n^2 \\
\end{bmatrix}z_x,
\end{equation}
and
\begin{equation}\label{eq:p3}
p_3(y,z_y)=n^2\gamma z_y^T\begin{bmatrix}
5y_1^2 & 2y_1y_2 & \cdots & 2y_1y_n  \\
2y_1y_2 & 5y_2^2  & \cdots & 2y_2y_n \\
\vdots & \aaan{\vdots} &\ddots &\vdots \\
2y_1y_n &  \cdots &2y_{n-1}y_n & 5y_n^2 \\
\end{bmatrix}z_y.
\end{equation}

We show that (\ref{eq:p1+p2+p3}) is sos by showing that $p_1$,
$p_2$, and $p_3$ are each individually sos. %
\aaa{To see that $p_2$ is sos, simply note that we can rewrite it
as \alex{\[ p_2(x,z_x) = n^2 \gamma \left[ 3 \sum_{k=1}^n
z_{x,k}^2 x_k^2 + 2 \Big(\sum_{k=1}^n z_{x,k} x_k\Big)^2 \right].
\]} The argument for $p_3$ is of course identical.} \alex{To show that $p_1$ is sos,} we argue as follows. If we
multiply out \aaan{the first term} $2z_x^TC(x,y)z_y$, we
\jnt{obtain} a polynomial \alex{with monomials of the form}
\begin{equation}\label{eq:typical.monomial}
\pm2\beta_{i,j,k,l}z_{x,k}x_iy_jz_{y,l},
\end{equation}
where $0\leq \beta_{i,j,k,l} \leq \gamma$, \alex{by the definition
of} $\gamma$. Since
\begin{equation}
\pm2\beta_{i,j,k,l}z_{x,k}x_iy_jz_{y,l}+\beta_{i,j,k,l}z_{x,k}^2x_i^2+\beta_{i,j,k,l}y_j^2z_{y,l}^2=\beta_{i,j,k,l}(z_{x,k}x_i\pm
y_jz_{y,l})^2,
\end{equation} \alex{by pairing up the terms of $2z_x^T C(x,y) z_y$
with \aaa{fractions of} the squared terms $z_{x,k}^2 x_i^2$ and
$z_{y,l}^2 y_j^2$, we get a sum of squares. Observe that there are
more than enough squares for each monomial of $2z_x^T C(x,y) z_y$
because each such monomial $\pm2\beta_{i,j,k,l}z_{x,k}x_i y_j
z_{y,l}$ occurs at most once, so that each of the terms $z_{x,k}^2
x_i^2$ and $z_{y,l}^2 y_j^2$ will be needed at most $n^2$ times,
each time with a coefficient of at most $\gamma$.} Therefore,
$p_1$ is sos, \alex{and this completes the proof.}
\end{proof}
\aaa{We can now \aaan{complete the proof of} strong NP-hardness of
deciding convexity of quartic forms.
\begin{proof}[Proof of Theorem~\ref{thm:convexity.quartic.nphard}]
As we remarked earlier, deciding nonnegativity of biquadratic
forms is known to be strongly
NP-hard~\cite{Ling_et_al_Biquadratic}. Given such a biquadratic
form $b(x;y)$, we can construct the polynomial $f(x,y)$ as in
(\ref{eq:f(x,y).defn.}). Note that $f(x,y)$ has degree four and is
homogeneous. Moreover, the reduction from $b(x;y)$ to $f(x,y)$
\aaan{runs in polynomial time} as we are only adding to $b(x;y)$
$2n+2 \binom{n}{2}$ new monomials with coefficient
$\frac{n^2\gamma}{2}$, and the size of $\gamma$ is by definition
only polynomially larger than the size of any coefficient of
$b(x;y)$. Since by Theorem~\ref{thm:main.reduction} convexity of
$f(x,y)$ is equivalent to nonnegativity of $b(x;y)$, we conclude
that deciding convexity of quartic forms is strongly NP-hard.
\end{proof}}
\paragraph {An algebraic version of the reduction.} Before we
proceed further with our results, we \jnt{make} a slight detour
and present an algebraic \aaan{analogue} of this reduction,
\jnt{which relates} sum of squares biquadratic forms to sos-convex
polynomials. Both of these concepts are well-studied in the
literature, in particular in regards to their connection to
semidefinite programming; see,
e.g.,~\cite{Ling_et_al_Biquadratic},~\cite{AAA_PP_not_sos_convex_journal},
and references therein.

\begin{definition}\label{defn:sos-convex}
A polynomial $p(x)$, with its Hessian denoted by $H(x)$, is
\emph{sos-convex} if the polynomial $y^TH(x)y$ is a sum of squares
in variables (x;y).\footnote{See~\cite{AAA_PP_CDC10_algeb_convex}
for three other equivalent definitions of sos-convexity.}
\end{definition}

\begin{theorem}\label{thm:algeb.ver.of.reduction}
Given a biquadratic form \aaa{$b(x;y)$}, let $f(x,y)$ be the
quartic form defined as in (\ref{eq:f(x,y).defn.}). Then $b(x;y)$
is a sum of squares if and only if $f(x,y)$ is sos-convex.
\end{theorem}
\begin{proof}
The proof is very similar to the proof of
\aaa{Theorem~\ref{thm:main.reduction}} and is left to the reader.
\end{proof}

Perhaps of independent interest,
\aaa{Theorems~\ref{thm:main.reduction}} and
\ref{thm:algeb.ver.of.reduction} imply that our reduction gives an
explicit way of constructing convex but not sos-convex quartic
forms (see~\cite{AAA_PP_not_sos_convex_journal}), starting from
any example of a psd but not sos biquadratic form
(see~\cite{Choi_Biquadratic}).

\paragraph {Some NP-hardness \jnt{results, obtained} as corollaries.} NP-hardness of checking convexity of quartic
forms directly establishes NP-hardness\aaa{\footnote{\aaa{All of
our NP-hardness results in this paper are in the strong sense. For
the sake of brevity, from now on we refer to strongly NP-hard
problems simply as NP-hard problems.}}}
 of \jnt{several problems of interest.}
Here, we mention a \jnt{few examples.}

\begin{corollary}\label{cor:nonnegativity.quartic.nphard}
It is NP-hard to decide \alex{nonnegativity of a homogeneous}
polynomial $q$ of degree four, \alex{of the form}
\begin{equation}\nonumber
q(x,y)=\frac{1}{2}p(x)+\frac{1}{2}p(y)-p\left(\textstyle{\frac{x+y}{2}}\right),
\end{equation}
for some homogeneous quartic polynomial $p$.
\end{corollary}
\begin{proof}
\jnt{Nonnegativity of $q$ is equivalent to convexity of $p$, and the result}
follows directly from
Theorem~\ref{thm:convexity.quartic.nphard}.%
\end{proof}

\begin{definition}
A set $\mathcal{S}\subset\mathbb{R}^n$ is \emph{basic closed
semialgebraic} if it can be written as
\begin{equation}\label{eq:semialgebraic.set}
\mathcal{S}=\{x\in\mathbb{R}^n|\ f_i(x)\geq0,\ i=1,\ldots,m\},
\end{equation}
for some positive integer $m$ and some polynomials $f_i(x)$.
\end{definition}

\begin{corollary}\label{cor:semialgeb.set.nphard}
Given a basic closed semialgebraic set $\mathcal{S}$ as in
(\ref{eq:semialgebraic.set}), where at least one of the defining
polynomials $f_i(x)$ has degree four, it is NP-hard to decide whether
$\mathcal{S}$ is a convex set.
\end{corollary}
\begin{proof}
Given a quartic polynomial $p(x)$, consider the basic closed
semialgebraic set $$\mathcal{E}_p=\{(x,t)\in\mathbb{R}^{n+1}|\
t-p(x)\geq0\},$$ describing the epigraph of $p(x)$. Since $p(x)$
is convex if and only if its epigraph is a convex set, the result
follows.\footnote{Another proof of this corollary is given by the
NP-hardness of checking convexity of sublevel sets of quartic
polynomials (Theorem~\ref{thm:quasi.pseudo.quartic.nphard} in
Section~\ref{subsec:quasi.pseudo.hard.degrees}).}
\end{proof}

\paragraph {Convexity of polynomials of even degree larger than
four.} We end this section by extending our hardness result to
polynomials of higher degree.
\begin{corollary}
It is NP-hard to check convexity of polynomials of any fixed even
degree $d\geq4$.
\end{corollary}
\begin{proof}
We have already established the result for polynomials of degree
four. Given such a degree four polynomial
$p(x)\mathrel{\mathop:}=p(x_1,\ldots,x_n)$ and an even degree
$d\geq6$, consider the polynomial $$q(x,x_{n+1})=p(x)+x_{n+1}^d$$
in $n+1$ variables. It is clear (e.g., from the block diagonal
structure of the Hessian of $q$) that $p(x)$ is convex if and only
if $q(x)$ is convex. The \jnt{result} follows.
\end{proof}

\section{Complexity of deciding strict convexity and strong
convexity}\label{sec:strict.strong}

\subsection{Definitions and basics}\label{subsec:strict.strong.basics}

\begin{definition}
A function $f:\mathbb{R}^n\rightarrow\mathbb{R}$ is \emph{strictly
convex} if %
for all $x\neq y$ %
 and all $\lambda \in (0,1)$, we have
\begin{equation}\label{eq:strict.convexity.defn.}
f(\lambda x+(1-\lambda)y)< \lambda f(x)+(1-\lambda)f(y).
\end{equation}
\end{definition}

\begin{definition}
A twice differentiable function
$f:\mathbb{R}^n\rightarrow\mathbb{R}$ is \emph{strongly convex} if
its Hessian $H(x)$ satisfies
\begin{equation}\label{eq:strong.convexity.defn}
H(x)\succeq mI,
\end{equation}
for a scalar $m>0$ and for all $x$.%
\end{definition}

We have the standard implications
\begin{equation}\label{eq:implications.strong.strict.convex}
\mbox{strong convexity} \ \Longrightarrow \  \mbox{strict
convexity}\ \Longrightarrow \ \mbox{convexity},
\end{equation}
but \jnt{none of the converse implications is true.}

\subsection{Degrees that are easy}\label{subsec:strict.strong.easy.degrees}
From the implications in
(\ref{eq:implications.strong.strict.convex}) and our previous
discussion, it is clear that odd degree polynomials can never be
strictly convex \aaan{or} strongly convex. \aaa{We cover the case
of quadratic polynomials in the following straightforward
proposition.}

\aaa{
\begin{proposition}%
For a quadratic polynomial $p(x)=\frac{1}{2}x^TQx+q^Tx+c$, the
notions of strict convexity and strong convexity are equivalent,
and can be decided in polynomial time.
\end{proposition}
}
\begin{proof}
Strong convexity always implies strict convexity. For the
\jnt{reverse} direction, assume that $p(x)$ is not strongly
convex. In view of (\ref{eq:strong.convexity.defn}), this means
that the matrix $Q$ is not positive definite. If $Q$ has a
negative eigenvalue, $p(x)$ is not convex, let alone strictly
convex. If $Q$ has a zero eigenvalue, let $\bar{x}\neq 0$ be the
corresponding eigenvector. Then $p(x)$ restricted to the
\aaan{line} from the origin to $\bar{x}$ is linear and hence not
strictly convex.

To see that these properties can be checked in polynomial time,
note that $p(x)$ is strongly convex if and only if the symmetric
matrix $Q$ is positive definite. \aaan{By Sylvester's criterion,}
positive definiteness of an $n \times n$ symmetric matrix is
equivalent to positivity of its $n$ leading principal minors,
\alex{each of which can be computed in polynomial time.}
\end{proof}

\subsection{Degrees that are hard}\label{subsec:strict.strong.hard.degrees}
With little effort, we can extend our NP-hardness result in the
previous section to \jnt{address} strict convexity and strong
convexity.

\begin{proposition}\label{thm:strong.convexity.NPhard.even.deg} %
It is NP-hard to decide strong convexity of polynomials of any
fixed even degree $d\geq4$.
\end{proposition}

\begin{proof}
\aaa{We give a reduction from the problem of deciding convexity of
quartic forms. Given a homogenous quartic polynomial
$p(x)\mathrel{\mathop:}=p(x_1,\ldots,x_n)$ and an even degree
$d\geq4$, consider the polynomial
\begin{equation}\label{eq:q.in.reduction.strong.convexity}
q(x,x_{n+1}):=p(x)+x_{n+1}^d+\textstyle{\frac{1}{2}}(x_1^2+\cdots+x_n^2+x_{n+1}^2)
\end{equation}
in $n+1$ variables. We claim that $p$ is convex if and only if $q$
is strongly convex. Indeed, if $p(x)$ is convex, then so is
$p(x)+x_{n+1}^d$. Therefore, the Hessian of $p(x)+x_{n+1}^d$ is
PSD. On the other hand, the Hessian of the term
$\frac{1}{2}(x_1^2+\cdots+x_n^2+x_{n+1}^2)$ is the identity
matrix. So, the minimum eigenvalue of the Hessian of
$q(x,x_{n+1})$ is positive and bounded \aaan{below by one}. Hence,
$q$ is strongly convex.

Now suppose that $p(x)$ is not convex. Let us denote the Hessians
of $p$ and $q$ respectively by $H_p$ and $H_q$. If $p$ is not
convex, then there exists a point $\bar{x}\in\mathbb{R}^n$ such
that $$\lambda_{\min}(H_p(\bar{x}))<0,$$ where $\lambda_{\min}$
here denotes the minimum eigenvalue. Because $p(x)$ is homogenous
of degree four, we have
$$\lambda_{\min}(H_p(c\bar{x}))=c^2\lambda_{\min}(H_p(\bar{x})),$$
for any scalar $c\in\mathbb{R}$. Pick $c$ large enough such that
$\lambda_{\min}(H_p(c\bar{x}))<1$. Then it is easy to see that
$H_q(c\bar{x},0)$ has a negative eigenvalue and hence $q$ is not
convex, let alone strongly convex.}
\end{proof}

\begin{remark}
It is worth noting that homogeneous polynomials of degree $d>2$
can never be strongly convex (because their Hessians vanish at the
origin). Not surprisingly, the polynomial $q$ in the proof of
\aaa{Proposition}~\ref{thm:strong.convexity.NPhard.even.deg} is
not homogeneous.
\end{remark}

\begin{proposition}
It is NP-hard to decide strict convexity of polynomials of any
fixed even degree $d\geq4$.
\end{proposition}

\begin{proof}
The proof is almost identical to the proof of
\aaa{Proposition}~\ref{thm:strong.convexity.NPhard.even.deg}. Let
$q$ be defined as in (\ref{eq:q.in.reduction.strong.convexity}).
If $p$ is convex, then we established that $q$ is strongly convex
and hence also strictly convex. If $p$ is not convex, we showed
that $q$ is not convex and hence also not strictly convex.
\end{proof}

\section{Complexity of deciding quasiconvexity and
pseudoconvexity}\label{sec:quasi.pseudo}
\subsection{Definitions and basics}\label{subsec:quasi.pseudo.basics}

\begin{definition}\label{def:quasiconvex.fn}
A function $f:\mathbb{R}^n\rightarrow\mathbb{R}$ is
\emph{quasiconvex} if \aaan{its sublevel sets}
\begin{equation}\label{eq:sublevel.sets}
\mathcal{S}(\alpha)\mathrel{\mathop:}=\{x\in\mathbb{R}^n \mid
f(x)\leq\alpha\},
\end{equation}for all $\alpha\in\mathbb{R}$, \aaan{are} convex.
\end{definition}

\begin{definition}\label{def:pseudoconvex.fn}
A differentiable function $f:\mathbb{R}^n\rightarrow\mathbb{R}$ is
\emph{pseudoconvex} if the implication
\begin{equation}\label{eq:pseudoconvexity.defn}
\nabla f(x)^T(y-x)\geq0  \ \Longrightarrow \ f(y)\geq f(x)
\end{equation}
holds for all $x$ and $y$ \aaa{in $\mathbb{R}^n$}.
\end{definition}

The following implications are well-known \aaa{(see e.g.~\cite[p.
143]{NLP_Book_Bazaraa})}:
\begin{equation}\label{eq:implications.convex.quasi.pseudo}
\mbox{convexity} \ \Longrightarrow \  \mbox{pseudoconvexity}\
\Longrightarrow \ \mbox{quasiconvexity},
\end{equation}
but the converse of neither implication is true in general.

\subsection{Degrees that are easy}\label{subsec:quasi.pseudo.easy.degrees}
\aaa{As we remarked earlier, linear polynomials are always convex
and hence also pseudoconvex and quasiconvex. Unlike convexity,
however, it is possible for polynomials of odd degree $d\geq3$ to
be pseudoconvex or quasiconvex.}
We will show in this section that \alex{somewhat} surprisingly,
quasiconvexity and pseudoconvexity of polynomials of any fixed odd
degree can be decided in polynomial time. Before we present these
results, we will \alex{cover} the easy case of quadratic
polynomials.

\aaa{
\begin{proposition}\label{thm:conv=quasi=pseudo.for.quadratics}
For a quadratic polynomial $p(x)=\frac{1}{2}x^TQx+q^Tx+c$, the
notions of convexity, pseudoconvexity, and quasiconvexity are
equivalent, and can be decided in polynomial time.
\end{proposition}
}
\begin{proof}
We \alex{argue} that the quadratic polynomial \alex{$p(x)$} is
convex if and only if
it is quasiconvex. %
\alex{Indeed, if $p(x)$ is not convex, then $Q$ has a negative
eigenvalue; letting $\bar{x}$ be \aaan{a} corresponding
eigenvector, we have that $p(t \bar{x})$ is a \jnt{quadratic}
polynomial in $t$, with negative leading coefficient, so $p(t
\bar{x})$ is not quasiconvex, \jnt{as a function of $t$.} This,
however, implies that $p(x)$ is not quasiconvex.}

\aaa{We have already argued in
Section~\ref{subsec:convexity.easy.degrees} that convexity of
quadratic polynomials can be decided in polynomial time.}
\end{proof}

\subsubsection{Quasiconvexity of polynomials of odd degree}
\label{subsubsec:quasi.odd.degree} In this subsection, we provide
a polynomial time algorithm for checking whether an odd-degree
polynomial is quasiconvex. \aaan{Towards this goal, we will first
show that quasiconvex polynomials of odd degree have a very
particular structure
(Proposition~\ref{prop:pseudo.odd.degree.special.struc.}).}

Our first lemma concerns quasiconvex polynomials of odd degree in
\aaan{one} variable. \jntn{A version of this lemma is provided
in~\cite[p.\ 99]{BoydBook}, without proof. We provide here a
proof, for completeness.}

\begin{lemma} \label{lemma:quasiconvex-odd-degree} Suppose that
$p(t)$ is a quasiconvex univariate polynomial of odd degree. Then,
$p(t)$ is monotonic.
\end{lemma}
\begin{proof} Suppose that monotonicity does not hold for the points $a < b < c$. This can happen
in two ways: \[ p(a)<p(b)\ \ {\rm and}\ \ p(c)<p(b),\] or \[ p(a)
> p(b)\ \ {\rm and}\ \   p(c)>p(b).
\] In the first case, pick $\alpha$ such that $p(b) > \alpha >
\max(p(a),p(c))$. The \aaan{sublevel} set $\{x\mid p(x) \leq
\alpha \}$ includes $a$ and $c$, but not $b$, and is therefore
non-convex, which contradicts the quasiconvexity of $p$.

In the second case, pick $\alpha$ such that $\min(p(a),p(c)) >
\alpha > p(b)$, and consider the \aaan{sublevel} set
$\aaan{\mathcal{S}(\alpha)} = \{ x \mid p(x) \leq \alpha \}$. If
the leading coefficient of $p$ is positive, then
$\aaan{\mathcal{S}(\alpha)}$ contains contains $b$ and some point
to the left of $a$ (since $\lim_{t \rightarrow -\infty} p(t) =
-\infty$) but not $a$, which is a contradiction. Similarly, if the
leading coefficient of $p$ is negative, then
$\aaan{\mathcal{S}(\alpha)}$ contains $b$ and some point to the
right of $c$ but not $c$, which is a contradiction.
\end{proof}

Next, we use the \aaan{preceding lemma} to characterize the
complements of \aaan{sublevel} sets of quasiconvex polynomials of
odd degree.

\begin{lemma} \label{lemma:halfspaces} Suppose that $p(x)$
is a quasiconvex polynomial of odd degree $d$. Then the set $\{ x
\mid  p(x) \geq \alpha \}$ is convex.
\end{lemma}

\begin{proof} Suppose not. In that case, there exist $x,y,z$ such that
$z$ is on the line segment connecting $x$ and $y$, \jntn{and such
that} $p(x), p(y) \geq \alpha$ but $p(z) < \alpha$. Consider the
polynomial
\[ q(t)=p(x + t(y-x)).\] This is, of course, a quasiconvex
polynomial with $q(0)=p(x)$, $q(1)=p(y)$, and $q(t')=p(z)$, for
some $t' \in (0,1)$. If $q(t)$ has degree $d$, then, by Lemma
\ref{lemma:quasiconvex-odd-degree}, it must be monotonic, which
immediately provides a contradiction.

Suppose now that $q(t)$ has degree less than $d$. Let us attempt
to perturb $x$ to $x+x'$, and $y$ to $y+y'$, so that the new
polynomial \[ \hat{q}(t)=p \left( x+x' + t(y+y'-x-x') \right)\]
has the following two properties: (i) $\hat{q}(t)$ is a polynomial
of degree $d$, and (ii) $\hat{q}(0)
> \hat{q}(t')$, $\hat{q}(1) > \hat{q}(t')$. If such perturbation
vectors $x', y'$ can be found, then we obtain a contradiction as
in the previous paragraph.

To satisfy condition (ii), it suffices (by continuity) to take
$x',y'$ with $\|x'\|, \|y'\|$ small enough. Thus, we \jntn{only}
need to argue that we can find arbitrarily small $x',y'$ that
satisfy condition (i). Observe that  the \jntn{coefficient of
$t^d$ in the polynomial} $\hat{q}(t)$ is a nonzero polynomial in
$x+x', y+y'$; let us denote \jntn{that coefficient} as
$\jntn{r}(x+x',y+y')$. \jntn{Since $r$ is a nonzero polynomial, it
cannot vanish at all points of any given ball. Therefore, even
when considering a small ball around $(x,y)$ (to satisfy condition
(ii)), we can find $(x+x',y+y')$ in that ball, with
$r(x+x',y+y')\neq 0$, thus establishing that the degree of $\hat
q$ is indeed $d$.} This completes the proof.

\old{Since the only polynomial which is zero on a ball of positive
radius is the zero polynomial, and $l$ is nonzero, we have that no
matter how small a ball ${\cal B}$ \aaa{we} take around $(x,y)$
there are some $x',y'$ with $(x+x',y+y') \in {\cal B}$ and
$l(x+x', y+y') \neq 0$. \aaa{[* Can you rephrase this argument?
Right now it is a bit confusing because it says that we know the
leading coefficient is non-zero but then for any ball we can take
the leading coefficient to be non-zero .... ?? *]} This completes
the proof.}
\end{proof}

\jntn{We now proceed to \aaan{a characterization} of quasiconvex
polynomials of odd degree.}

\begin{proposition}\label{prop:pseudo.odd.degree.special.struc.} %
Let $p(x)$ be a polynomial of odd degree $d$. Then, $p(x)$ is
quasiconvex if and only if it can be written as
\begin{equation} \label{quasiconvexdecomposition}
p(x)=h(\xi^Tx),
\end{equation}
for some \jntn{nonzero} $\xi \in \mathbb{R}^n$, and for some
\jntn{monotonic univariate} polynomial \jntn{$h(t)$ of degree $d$.
If, in addition, we require the nonzero component of $\xi$ with
the smallest index to be equal to unity, then $\xi$ and $h(t)$ are
uniquely determined by $p(x)$.}
\end{proposition}

\begin{proof} \old{\jntn{PROOF REWRITTEN}} It is easy to see that any polynomial
that can be written in the above form is quasiconvex. In order to
prove the converse, let us assume that $p(x)$ is quasiconvex. By
the definition of quasiconvexity, the closed set
$\aaan{\mathcal{S}}(\alpha)=\{x\mid p(x) \leq \alpha \}$ is
convex. On the other hand, Lemma \ref{lemma:halfspaces} states
that the closure of the complement of $\aaan{\mathcal{S}}(\alpha)$
is also convex. It is not hard to verify that, as a consequence of
these two properties, the set $\aaan{\mathcal{S}}(\alpha)$ must be
a halfspace. Thus, for any given $\alpha$, the sublevel set
$\aaan{\mathcal{S}}(\alpha)$ can be written as $\{x\mid
\xi(\alpha)^T x \leq c\aaan{(\alpha)}\}$ for some $\xi(\alpha) \in
\mathbb{R}^n$ and $ c(\alpha) \in \mathbb{R}$. This of course
implies that the level sets $\{x\mid p(x) = \alpha \}$ are
hyperplanes of the form $\{x\mid \xi(\alpha)^T x = c(\alpha)\}$.

We note that the sublevel sets are necessarily nested: if
$\alpha<\beta$, then $\aaan{\mathcal{S}}(\alpha)\subseteq
\aaan{\mathcal{S}}(\beta)$. An elementary consequence of this
property is that the hyperplanes must be collinear, i.e., that the
vectors $\xi(\alpha)$ must be positive multiples of each other.
Thus, by suitably scaling the coefficients $c(\alpha)$, we can
assume, without loss of generality, that $\xi(\alpha)=\xi$, for
some $\xi \in \mathbb{R}^n$, and for all $\alpha$. We then have
that $\{x\mid p(x)=\alpha\}=\{x\mid \xi^T x= c(\alpha)\}$.
Clearly, there is a one-to-one correspondence between $\alpha$ and
$c(\alpha)$, and therefore the value of $p(x)$ is completely
determined by $\xi^T x$. In particular, there exists a function
$h(t)$ such that $p(x)=h(q^Tx)$. Since $p(x)$ is a polynomial of
degree $d$, it follows that $h(t)$ is a univariate polynomial of
degree $d$. Finally, we observe that if $h(t)$ is not monotonic,
then $p(x)$ is not quasiconvex. This proves that a representation
of the desired form exists. Note that by suitably scaling $\xi$,
we can also impose the condition that the nonzero component of
$\xi$ with the smallest index is equal to one.

Suppose that now that $p(x)$ can also be represented in the form
$p(x)=\bar h(\bar \xi^T x)$ for some other polynomial $\bar h(t)$
and vector $\aaan{\bar{\xi}}$. Then, the gradient vector of $p(x)$
must be proportional to both $\xi$ and $\bar\xi$. The vectors
$\xi$ and $\bar \xi$ are therefore collinear. Once we impose the
requirement that the nonzero component of $\xi$ with the smallest
index is equal to one, we obtain that $\xi=\bar\xi$ and,
consequently, $h=\bar h$. This establishes the claimed uniqueness
of the representation.
\end{proof}

\old{ \alex{\noindent {\it \aaa{Remark.}} Observe that if Eq.\
(\ref{quasiconvexdecomposition}) holds for one \jnt{choice of}
$q,h$, then it holds for infinitely many \jnt{choices of $q,h$.}
For example, one can multiply every entry of the vector $q$ by
some number $c$ while multiplying the $i$th coefficient of $h$ by
$1/c^i$. If the original $h$ is monotonic, the new $h$ will be as
well, since it can be written as $h(x/c)$. In particular, we may
assume without loss of generality that $h$ is monic.
Alternatively, we may assume without loss of generality that
$q_1=1$; this may be assured by the above transformation unless it
happens that $q_1=0$, in which case $p(x)$ does not depend on
$x_1$.} }

\bigskip

\noindent {\it Remark.} It is not hard to see that if $p(x)$ is
homogeneous and quasiconvex, then one can additionally conclude
that $h(t)$ can be taken to be $h(t)=t^d$, where $d$ is the degree
of $p(x)$.

\begin{theorem}\label{thm:quasi.odd.degree.poly.time}
\jntn{For any fixed odd degree $d$, the} quasiconvexity of
polynomials of degree $d$ can be checked in polynomial time.
\end{theorem}

\begin{proof} \old{\jntn{NEW PROOF}}
The algorithm consists of attempting to build a representation of
$p(x)$ of the form given in Proposition
\ref{prop:pseudo.odd.degree.special.struc.}. The polynomial $p(x)$
is quasiconvex if and only if the attempt is successful.

Let us proceed under the assumption that $p(x)$ is quasiconvex. We
differentiate $p(x)$ symbolically to obtain its gradient vector.
Since a representation of the form given in Proposition
\ref{prop:pseudo.odd.degree.special.struc.} exists, the gradient
is of the form $\nabla p(x) =\xi h'(\xi^T x)$, where $h'(t)$ is
the derivative of $h(t)$. In particular, the different components
of the gradient are polynomials that are proportional to each
other. (If they are not proportional, we conclude that $p(x)$ is
not quasiconvex, and the algorithm terminates.) By considering the
ratios between different components, we can identify the vector
$\xi$, up to a scaling factor. By imposing the additional
requirement that the nonzero component of $\xi$ with the smallest
index is equal to one, we can identify $\xi$ uniquely.

We now proceed to identify the polynomial $h(t)$. For
$k=1,\ldots,d+1$, we evaluate $p(k\aaan{\xi})$, which must be
equal to $h(\xi^T\xi  k)$. We thus obtain the values of $h(t)$ at
$d+1$ distinct points, from which $h(t)$ is completely determined.
We then verify that $h(\xi^T x)$ is indeed equal to $p(x)$. This
is easily done, in polynomial time, by writing out the $O(n^d)$
coefficients of these two polynomials in $x$ and verifying that
they are equal. (If they are not all equal, we conclude that
$p(x)$ is not quasiconvex, and the algorithm terminates.)

Finally, we test whether the above constructed univariate
polynomial $h$ is monotonic, i.e., whether its derivative $h'(t)$
is either nonnegative or nonpositive. This can be accomplished,
e.g., by quantifier elimination or by other well-known algebraic
techniques for counting the number and the multiplicity of real
roots of univariate polynomials;
see~\cite{Algo_real_algeb_geom_Book}. Note that this requires only
a constant number of arithmetic operations since the degree $d$ is
fixed. If $h$ fails this test, then $p(x)$ is not quasiconvex.
Otherwise, our attempt has been successful and we decide that
$p(x)$ is indeed quasiconvex.
\end{proof}

\subsubsection{Pseudoconvexity of polynomials of odd degree}
\aaan{In analogy to
Proposition~\ref{prop:pseudo.odd.degree.special.struc.}, we
present next a characterization of odd degree pseudoconvex
polynomials, which gives rise to a polynomial time algorithm for
checking this property.}

\alex{\begin{corollary}\label{thm:quasi.odd.degree.special.struc.} %
Let $p(x)$ be a polynomial of odd degree $d$. Then, $p(x)$ is
pseudoconvex if and only if $p(x)$ can be written in the form
\begin{equation} \label{pseudoconvexrepresentation}
p(x)=h(\aaan{\xi}^T x),
\end{equation}
for some $\aaan{\xi} \in \mathbb{R}^n$ and some univariate
polynomial $h$ of degree $d$ such that \aaa{its derivative}
$h'(t)$ has no real roots.
\end{corollary}}

\alex{ \noindent {\it \aaa{Remark.}} Observe that polynomials $h$
with $h'$ having no real roots \jnt{comprise a subset of the set}
of monotonic polynomials.}

\begin{proof} \alex{\jnt{Suppose that $p(x)$ is pseudoconvex.}
Since a pseudoconvex polynomial is quasiconvex, it admits a
representation $h(\aaan{\xi}^T x)$ where $h$ is monotonic. If
$h'(t)=0$ for some $t$, then picking $a = t \cdot
\aaan{\xi}/\|\aaan{\xi}\|_2^2$, we have that $\nabla p(a)=0$, so
that by pseudoconvexity, $p(x)$ is minimized at $a$. This,
however, is impossible \aaa{since an odd degree polynomial is
never bounded below}.}
\alex{Conversely, suppose $p(x)$ can be represented as in Eq.
(\ref{pseudoconvexrepresentation}). \jnt{Fix some} $x,y$, and
define the polynomial $\aaa{u}(t)=p(x+t(y-x))$. Since
$\aaa{u}(t)=h(\aaan{\xi}^T x + t \aaan{\xi}^T (y-x))$, we have
that either (i) $\aaa{u}(t)$ is constant, or (ii) $\aaa{u}'(t)$
has no real roots. Now if $\nabla p(x)(y-x) \geq 0$, then
$\aaa{u}'(0) \geq 0$. Regardless of whether (i) or (ii) holds,
this implies that $\aaa{u}'(t) \geq 0$ everywhere, so that
$\aaa{u}(1) \geq \aaa{u}(0)$ or $p(y) \geq p(x)$.}
 \end{proof}
\begin{corollary}\label{thm:pseudo.odd.degree.poly.time}
\aaan{For any fixed odd degree $d$, the} pseudoconvexity of
polynomials of degree $d$ can be checked in polynomial time.
\end{corollary}
\begin{proof} \alex{This is a simple modification of our algorithm
for testing quasiconvexity
(Theorem~\ref{thm:quasi.odd.degree.poly.time}). \aaan{The first
step of the algorithm is in fact identical: once we impose the
additional requirement that the nonzero component of $\xi$ with
the smallest index should be equal to one, we can uniquely
determine the vector $\xi$ and the coefficients of the univariate
polynomial $h(t)$ that satisfy Eq.
(\ref{pseudoconvexrepresentation}) . (If we fail, $p(x)$ is not
quasiconvex and hence also not pseudoconvex.) Once we have $h(t)$,
we can check whether $h'(t)$ has no real roots e.g. by computing
the signature of the Hermite \aaa{form of $h'(t)$;
see~\cite{Algo_real_algeb_geom_Book}}.}}

\end{proof}

\begin{remark}
Homogeneous polynomials of odd degree $d\geq3$ are never
pseudoconvex. The reason is that the gradient of these polynomials
vanishes at the origin, but yet the origin is not a global minimum
since odd degree polynomials are unbounded below.
\end{remark}

\subsection{Degrees that are hard}\label{subsec:quasi.pseudo.hard.degrees}
The main result of this section is the following theorem.

\begin{theorem}\label{thm:quasi.pseudo.quartic.nphard}
It is NP-hard to check quasiconvexity/pseudoconvexity of degree
four polynomials. This is true even when the polynomials are
restricted to be homogeneous.
\end{theorem}

In view of Theorem~\ref{thm:convexity.quartic.nphard}, which
established NP-hardness of deciding convexity of \emph{homogeneous}
quartic polynomials,
Theorem~\ref{thm:quasi.pseudo.quartic.nphard} follows immediately
from the following result.\footnote{A slight variant of
Theorem~\ref{thm:quasiconvexity.homog.same.convexity} has appeared
in~\cite{AAA_PP_CDC10_algeb_convex}.}

\begin{theorem}\label{thm:quasiconvexity.homog.same.convexity}
For a homogeneous polynomial $p(x)$ of even degree $d$, the notions
of convexity, pseudoconvexity, and quasiconvexity are all
equivalent.\footnote{The result is more generally true for
differentiable functions that are homogeneous of even degree. Also,
the requirements of homogeneity and having an even degree both
need to be present. Indeed, $x^3$ and $x^4-8x^3+18x^2$ are both
quasiconvex but not convex, the first being homogeneous of odd
degree and the second being nonhomogeneous of even degree.}
\end{theorem}

We start the proof of this theorem by first proving an easy lemma.

\begin{lemma} \label{lem:quasiconvex.homog.then.psd}
Let $p(x)$ be a quasiconvex homogeneous polynomial of even degree
\aaan{$d\geq2$}. Then $p(x)$ is nonnegative.
\end{lemma}

\begin{proof}
Suppose, to derive a contradiction, that there exist some
$\epsilon>0$ and $\bar{x}\in\mathbb{R}^n$ such that
$p(\bar{x})=-\epsilon$. Then by homogeneity of even degree we must
have $p(-\bar{x})=p(\bar{x})=-\epsilon$. On the other hand,
homogeneity of $p$ implies that $p(0)=0$. Since the origin is on
the line between $\bar{x}$ and $-\bar{x}$, this shows that the
\aaan{sublevel set} $\mathcal{S}(-\epsilon)$ is not convex,
contradicting the quasiconvexity of $p$.
\end{proof}

\begin{proof}[Proof of
Theorem~\ref{thm:quasiconvexity.homog.same.convexity}]
\alex{We show that a quasiconvex homogeneous polynomial of even degree is convex. In view of implication
(\ref{eq:implications.convex.quasi.pseudo}), this \alex{proves the theorem.}}

\alex{Suppose that $p(x)$ is a quasiconvex polynomial. Define
$\mathcal{S}=\{x\in\mathbb{R}^n \mid\  p(x)\leq1\}$. By
homogeneity, for any $a \in \mathbb{R}^n$ with $p(a)>0$, we have
that \[ \frac{a}{p(a)^{1/d}} \in \mathcal{S}.\] By quasiconvexity,
this \jnt{implies} that for any $a,b$ with $p(a),p(b)>0$, any
point on the line connecting $a/p(a)^{1/d}$ and $b/p(b)^{1/d}$ is
in $\mathcal{S}$. In particular, consider
\[ c=\frac{a+b}{p(a)^{1/d}+p(b)^{1/d}}.\]  Because $c$ can be written as
\[ %
c = \left( \frac{p(a)^{1/d}}{p(a)^{1/d} + p(b)^{1/d}}\right) \left(
\frac{a}{p(a)^{1/d}} \right) + \left( \frac{p(b)^{1/d}}{p(a)^{1/d} +
p(b)^{1/d}} \right) \left( \frac{b}{p(b)^{1/d}} \right),\] we have
that $c \in \mathcal{S}$, \jnt{i.e., $p(c)\leq 1$. By homogeneity, this  inequality} can be restated as
\[ p(a+b) \leq (p(a)^{1/d} + p(b)^{1/d})^d,\] and therefore
\begin{equation} p\Big(\frac{a+b}{2}\Big) \leq \left( \frac{p(a)^{1/d} + p(b)^{1/d}}{2}\right) ^d \leq \frac{p(a)+p(b)}{2}, \label{convexity-where-positive} \end{equation}
where the last inequality is due to the convexity of $x^d$.}

\alex{Finally, note that for any polynomial $p$, the set \jnt{$\{x\mid p(x)
\neq 0$\}} is dense in $\mathbb{R}^n$ (here we again appeal to the fact that the only polynomial that
is zero on a ball of positive radius is the zero polynomial); and since $p$ is nonnegative due
to Lemma \ref{lem:quasiconvex.homog.then.psd}, the set
\jnt{$\{x\mid p(x)
> 0$\}} is dense in $\mathbb{R}^n$.
\jnt{Using the
continuity of $p$, it follows that
Eq. (\ref{convexity-where-positive}) holds not only when $a,b$
satisfy $p(a),p(b)>0$, but for all $a$, $b$. Appealing to the continuity of $p$ again, we see that for all $a,b$,
$p(\lambda a + (1-\lambda) b) \leq \lambda p(a) + (1-\lambda) p(b)$, for all $\lambda \in [0,1]$. This establishes} that $p$ is convex.}

\end{proof}

\paragraph{Quasiconvexity/pseudoconvexity of polynomials of even degree larger than four.}

\begin{corollary}\label{cor:quasi.nphard.d>=4}
It is NP-hard to decide quasiconvexity of polynomials of any fixed
even degree $d\geq4$.
\end{corollary}

\begin{proof}\label{thm:quasi.nphard.even.degree}
We have already proved the result for $d=4$. To establish the result for
even degree $d\geq6$, recall that we have established NP-hardness
of deciding convexity of homogeneous quartic polynomials. Given
such a quartic form $p(x)\mathrel{\mathop:}=p(x_1,\ldots,x_n)$,
consider the polynomial
\begin{equation}\label{eq:q.in.reduction.quasi.even.degree}
q(x_1,\ldots, x_{n+1}) = p(x_1,\ldots,x_n) + x_{n+1}^d.
\end{equation}
We claim that $q$ is quasiconvex if and only if $p$ is convex.
Indeed, if $p$ is convex, then obviously so is $q$, and therefore
$q$ is quasiconvex. Conversely, if $p$ is not convex, then by
Theorem~\ref{thm:quasiconvexity.homog.same.convexity}, it is not
quasiconvex. So, there exist points $a,b,c\in\mathbb{R}^n$, with
$c$ on the line connecting $a$ and $b$, such that $p(a)\leq1$,
$p(b)\leq1$, but $p(c)>1$. Considering points $(a,0)$, $(b,0)$,
$(c,0)$, we see that $q$ is not quasiconvex. It follows that it is
NP-hard to decide quasiconvexity of polynomials of even degree
four or larger.
\end{proof}

\begin{corollary}
It is NP-hard to decide pseudoconvexity of polynomials of any
fixed even degree $d\geq4$.
\end{corollary}

\begin{proof}
The proof is almost identical to the proof of
\aaa{Corollary~\ref{cor:quasi.nphard.d>=4}}. \aaa{Let} $q$ be
defined as in (\ref{eq:q.in.reduction.quasi.even.degree}). If $p$
is convex, then $q$ is convex and hence also pseudoconvex. If $p$
is not convex, we showed that $q$ is not quasiconvex and hence
also not pseudoconvex.
\end{proof}

\section{Summary and conclusions}\label{sec:summary.conclusions}
We studied the computational complexity of testing convexity and
some of its variants, for polynomial functions. The notions that
we considered and the implications among them are summarized
below:

\vspace{3mm} strong convexity $\Longrightarrow$ strict convexity
$\Longrightarrow$ convexity $\Longrightarrow$ pseudoconvexity
$\Longrightarrow$ quasiconvexity. \vspace{3mm}

Our complexity results as a function of the degree of the
polynomial are listed in Table~1. %

\begin{table} 
\begin{center}
\begin{tabular}{l||c|c|c|c}
\textbf{property vs. degree} & 1 & 2 & odd $\geq3$ & even $\geq4$ \\
  \hline \hline
   \noindent\large{strong convexity} & no & P & no & strongly NP-hard \\
  \large{strict convexity} & no & P & no & strongly NP-hard  \\
  \large{convexity} & yes & P & no & strongly NP-hard \\
  \large{pseudoconvexity} & yes & P & P & strongly NP-hard \\
  \large{quasiconvexity} & yes & P & P & strongly NP-hard \\
\end{tabular}
\end{center}
\caption{Summary of our complexity results. A yes (no) entry means that the question is trivial for that particular entry because the answer is always yes (no) independent of the input. By \aaan{P}, we mean that the problem can be solved in polynomial time.}
\label{table:summary}
\end{table}

We gave polynomial time algorithms for checking pseudoconvexity
and quasiconvexity of odd degree polynomials that can be useful in
many applications. Our negative results, on the other hand, imply
(\jnt{under P$\neq$NP}) the impossibility of a polynomial time (or
even pseudo-polynomial time) algorithm for testing any of the
properties listed in Table~1 for polynomials of even degree four
or larger. \aaa{Although the implications of convexity are very
significant in optimization theory, our results suggest that
\aaan{unless additional structure is present,} ensuring the mere
presence of convexity is likely an intractable task.} It is therefore
natural to wonder whether there are other properties of
optimization problems that share some of the attractive
consequences of convexity, but are easier to check.

Of course, NP-hardness of a problem does not stop us from studying it,
but on the contrary, stresses the need for finding good approximation
algorithms that can deal with a large number of instances
efficiently. As an example, semidefinite programming based relaxations
relying on algebraic concepts such as sum of squares decomposition of
polynomials currently seem to be very promising techniques for
recognizing convexity of polynomials and basic semialgebraic sets. It
would be useful to identify special cases where these relaxations are
exact or give theoretical bounds on their performance guarantees.

\bibliographystyle{abbrv}
\bibliography{pablo_amirali}

\end{document}